\documentclass[12pt]{amsart}

\begin{document}

\theoremstyle{plain}
    \newtheorem{theorem}{Theorem}
    \newtheorem{lemma}[theorem]{Lemma}
\theoremstyle{definition}
    \newtheorem{definition}[theorem]{Definition}
\theoremstyle{remark}
    \newtheorem{remark}[theorem]{Remark}

\title[Strongly primitive species with potentials: aims and limitations]{Strongly primitive species with potentials: aims and limitations\\
{\tiny Oberwolfach talk, December 2013}}
\author{Daniel Labardini-Fragoso}
\address{Daniel Labardini-Fragoso\newline
Instituto de Matem\'aticas, UNAM, Mexico}
\email{labardini@im.unam.mx}

\begin{abstract}

\noindent This is an extended abstract of my talk at the Oberwolfach Workshop ``Cluster Algebras and Related Topics'' (December 8 - 14, 2013). It is based on a joint work with A. Zelevinsky (arXiv:1306.3495).
\end{abstract}

\maketitle


\noindent
This talk is based on \cite{DLF_LZ1}. We define mutations of strongly primitive species with potentials, and show that if $B$ is any $n\times n$ integer matrix admitting a skew-symmetrizer with pairwise coprime diagonal entries, then for any given sequence of mutations which one wants to perform there exists a species realization $A$ of $B$ which admits a potential $S$ such that the mutations of $(A,S)$ along the given sequence are compatible with matrix mutation.

An $n\times n$ integer matrix $B$ is said to be \emph{skew-symmetrizable} if there exist positive integers $d_1,\ldots,d_n$ such that $DB$ is skew-symmetric, where $D=\operatorname{diag}(d_1,\ldots,d_n)$. Such $D$ is said to be a \emph{skew-symmetrizer} of $B$. 

A \emph{weighted quiver} is a pair $(Q,\mathbf{d})$, where $Q=(Q_0,Q_1,t,h)$ is a loop-free quiver and $\mathbf{d}=(d_i)_{i\in Q_0}$ is a tuple of positive integers attached to the vertices of $Q$. We call $\mathbf{d}$ the \emph{weight tuple} of $(Q,\mathbf{d})$.

If $B$ is a skew-symmetrizable matrix with skew-symmetrizer $D$, we define a weighted quiver  $(Q_B,\mathbf{d})$ as follows:  The vertex set of $Q_B$ is $\{1,\ldots,n\}$, and $Q_B$ has exactly
$$
\frac{\operatorname{gcd}(d_i,d_j)b_{ij}}{d_j}
$$
arrows from $j$ to $i$ whenever $b_{ij}\geq 0$. The tuple $\mathbf{d}$ is defined to be the tuple of diagonal entries of $D$.

If $D$ (and hence $\mathbf{d}$) is kept fixed, the assignment $B\mapsto(Q_B,\mathbf{d})$ is easily seen to be a bijection between the set of $n\times n$ integer matrices that are skew-symmetrized by $D$ and the set of 2-acyclic weighted quivers with weight tuple $\mathbf{d}$. This means that there is a weighted quiver counterpart of the notion of matrix mutation:

\begin{definition}\label{DLF:def-weighted-quiver-mutation} Let $(Q,\mathbf{d})$ be a
weighted quiver. For $k\in Q_0$, define the \emph{mutation of $(Q,\mathbf{d})$ with respect to $k$} to be the weighted quiver $\mu_k(Q,\mathbf{d})$ with vertex set $Q_0$ and weight tuple $\mathbf{d}$, obtained from as the result of performing the following 3-step procedure:
\begin{itemize}
\item[(Step 1)] For each pair of arrows $a:j\to k$ and $b:k\to i$ of $Q$, introduce
$\frac{\gcd(d_i,d_j)d_k}{\gcd(d_i,d_k)\gcd(d_k,d_j)}$ ``composite'' arrows from $j$ to
$i$;
\item[(Step 2)] replace each $c\in Q_1$ incident to $k$ with an arrow $c^*$ going in the opposite direction;
\item[(Step 3)] choose a maximal collection of disjoint 2-cycles and remove them.
\end{itemize}
\end{definition}

Let $d$ be the least common multiple of the tuple $\mathbf{d}=(d_{i})_{i\in Q_0}$, $F$ be a finite field, and $E$ be the unique degree-$d$ field extension of $F$. For $i\in Q_0$, set $F_i$ to be the unique degree-$d_i$ field subextension of $E/F$, and for every pair of vertices $i,j\in Q_0$, set
$$
A_{ij}=\bigoplus_{a:j\to i} F_i\otimes_{F_i\cap F_j}F_j.
$$
Define
$$
R=\bigoplus_{i\in Q_0}F_i \ \ \ \text{and} \ \ \ A=\bigoplus_{i,j\in Q_0}A_{ij}.
$$
Then $R$ is a semisimple ring and $A$ is an $R$-$R$-bimodule. We say that $A$ is the \emph{species of $(Q,\mathbf{d})$} over $E/F$.

\begin{remark} The data $((F_i)_{i\in Q_0},(A_{ij})_{i,j\in Q_0},(A_{ij}^*)_{i,j\in Q_0})$ constitutes a species (or \emph{modulation}) of the valued quiver of $B$ in the sense of Dlab-Ringel \cite{DLF_DR}, hence our use of the term ``species".
\end{remark}

The complete tensor algebra of $A$ over $R$ is called the \emph{complete path algebra} of $A$ and denoted $R\langle\hspace{-0.05cm}\langle A\rangle\hspace{-0.05cm}\rangle$. Thus we have
$$
R\langle\hspace{-0.05cm}\langle A\rangle\hspace{-0.05cm}\rangle=\prod_{\ell=0}^{\infty}{A^\ell}
$$
as an $R$-$R$-bimodule, where $A^\ell$ denotes the $\ell$-fold tensor product $A\otimes_R\ldots\otimes_R A$.
If an element $S\in \prod_{\ell=1}^{\infty}{A^\ell}$ satisfies $S=\sum_{i\in Q_0}e_iSe_i$, where $e_i$ is the idempotent sitting in the $i^{\operatorname{th}}$ component of $R$, we say that $S$ is a \emph{potential} on $A$ and that $(A,S)$ is a \emph{species with potential}.

From now on we assume that $(Q,\mathbf{d})$ is \emph{strongly primitive}, that is, we suppose that $\operatorname{gcd}(d_i,d_j)=1$ for all $i\neq j$. We shall say that $A$, the species of $(Q,\mathbf{d})$ over $E/F$, is strongly primitive as well, and that $(A,S)$ is a \emph{strongly primitive species with potential} whenever $S$ is a potential on $A$. We also assume that the characteristic of $F$ is congruent to $1$ modulo $d$. This implies that there exists an element $v\in E$ such that the set $\mathcal{B}=\{1,v,v^2,\ldots,v^{d-1}\}$ is an \emph{eigenbasis} of $E/F$, that is, an $F$-basis of $E$ consisting of eigenvectors of all elements of the Galois group $\operatorname{Gal}(E/F)$. This eigenbasis is quite useful to obtain an explicit description of the elements of $R\langle\hspace{-0.05cm}\langle A\rangle\hspace{-0.05cm}\rangle$; indeed, setting $\mathcal{B}_i=\mathcal{B}\cap F_i$ for $i\in Q_0$ we have:

\begin{lemma} For every $\ell\geq 0$, the set $\{\omega_0 a_1\omega_1a_2\ldots\omega_{\ell-1} a_\ell\omega_\ell \ | \ t(a_m)=h(a_{m+1})$ for $m=1,\ldots,\ell-1$, $\omega_0\in\mathcal{B}_{h(a_1)}$ and $\omega_{m}\in\mathcal{B}_{t(a_m)}$ for $m=1,\ldots,\ell\}$, whose elements we call \emph{paths of length $\ell$}, is a basis of $A^\ell$ as an $F$-vector space. Consequently, every element of $R\langle\hspace{-0.05cm}\langle A\rangle\hspace{-0.05cm}\rangle$ has a unique expression as a possibly infinite $F$-linear combination of paths. In particular, every potential has a unique expression as a possibly infinite $F$-linear combination of cyclic paths of positive length.
\end{lemma}

Given a 2-acyclic strongly primitive weighted quiver $(Q,\mathbf{d})$ and a vertex $k\in Q_0$, for each pair of arrows $a:j\to k$ and $b:k\to i$ of $Q$, the composite arrows introduced in the first step of the weighted-quiver mutation with respect to $k$ are denoted with the symbols $[b\omega a]$, where $\omega$ runs in the set $\mathcal{B}_k$.

We can state now the mutation rule for strongly primitive species with potentials.

\begin{definition} Let $(Q,\mathbf{d})$ be a 2-acyclic strongly primitive weighted quiver, and let $A$ be its species over $E/F$. For $k\in Q_0$, let $\widetilde{\mu}_k(A)$ denote the species over $E/F$ of the weighted quiver obtained from $(Q,\mathbf{d})$ by applying only the first two steps of weighted-quiver mutation. If $S$ is a potential on $A$, we define a potential $\widetilde{\mu}_k(S)$ on $\widetilde{\mu}_k(A)$ according to the formula
$$
\widetilde{\mu}_k(S)=[S]+\sum_{\overset{a}{\rightarrow}k\overset{b}{\rightarrow}}\sum_{\omega\in\mathcal{B}_k}\omega^{-1}b^*[b\omega a]a^*.
$$
The species with potential which is the \emph{reduced part} of $(\widetilde{\mu}_k(A),\widetilde{\mu}_k(S))$ will be called the \emph{mutation} of $(A,S)$ with respect to $k$ and denoted $\mu_k(A,S)$.
\end{definition}

We refer the reader to \cite{DLF_LZ1} for the definition of reduced parts as well as for a proof of their existence. Note that the underlying species of $\mu_k(A,S)$ is again strongly primitive, although its underlying weighted quiver may have oriented 2-cycles.

\begin{definition} Given a finite sequence $(k_1,\ldots,k_m)$ of vertices of $Q$, we say that $(A,S)$ is $(k_m,\ldots,k_1)$-\emph{non-degenerate} if the quivers underlying the species with potentials $(A,S)$, $\mu_{k_1}(A,S)$, $\mu_{k_2}\mu_{k_1}(A,S)$, $\ldots$, $\mu_{k_m}\ldots\mu_{k_2}\mu_{k_1}(A,S)$, are 2-acyclic (hence well-defined).
\end{definition}

The following is the main result of the talk.

\begin{theorem}\cite{DLF_LZ1}\label{DLF:main-result} If $(Q,\mathbf{d})$ is a 2-acyclic strongly primitive weighted quiver, then for every finite sequence $(k_1,\ldots,k_m)$ of vertices of $Q$ there exists a finite-degree field extension $K/F$ which is linearly disjoint from $E/F$ and has the property that the species $A_{KE/K}$ of $(Q,\mathbf{d})$ over $KE/K$ admits a potential $S$ such that $(A_{KE/K},S)$ is $(k_1,\ldots,k_m)$-non-degenerate.
\end{theorem}

Besides Theorem \ref{DLF:main-result}, mutations of strongly primitive species with potentials share many other properties with the mutations of quivers with potentials of Derksen-Weyman-Zelevinsky. For example, they are well-defined up to right-equivalence and involutive up to right-equivalence, and they preserve Jacobi-finiteness. Following the spirit of \cite{DLF_DWZ}, the notion of mutation is further lifted in \cite{DLF_LZ1} to the representation-theoretic level.

The class of skew-symmetrizable matrices whose associated weighted quivers are strongly primitive includes several examples of matrices that do not admit \emph{global unfoldings} whatsoever (a global unfolding is an unfolding which is compatible with all possible sequences of mutations). So, the species framework in \cite{DLF_LZ1} provides a representation-theoretic approach to (the skew-symmetrizable matrices of) several cluster algebras where approaches via unfoldings do not work.

Now, how about matrices that admit global unfoldings?, and how about species with potentials for weighted quivers that are not strongly primitive? The answer to the first question is provided by the work \cite{DLF_Demonet} of L.~Demonet, who developed an approach to mutations via \emph{group species with potentials}. The framework of group species differs from the one we have presented above in that a group species attaches group algebras to the vertices of $Q$ rather than fields; the bimodules attached to the arrows are hence also different in nature from the bimodules $A_{ij}$ above. Group species with potentials provide a representation-theoretic approach that can be successfully applied to those matrices that admit global unfoldings through group actions. Unfortunately, for matrices that do not admit global unfoldings one needs a framework different from the group species framework.

Regarding the second question: if $\operatorname{gcd}(d_i,d_j)$ is not assumed to be equal to 1 for all $i\neq j$, one can easily construct an example of a weighted quiver $(Q,\mathbf{d})$ that has a 2-cycle $ab$ which is cyclically equivalent to $0$ in $R\langle\hspace{-0.05cm}\langle A\rangle\hspace{-0.05cm}\rangle$, and this implies that  none of the arrows $a,b$ belongs to the Jacobian ideal of any potential on the species $A$ of $(Q,\mathbf{d})$ over $E/F$. This means that none of $a,b$ can be deleted from $Q$. Ultimately, this yields an example of a 2-acyclic species (not strongly primitive) such that, no matter which potential on it we take,  when we try to perform the three steps of weighted-quiver mutation at the level of species with potentials, the species framework presented above fails to delete 2-cycles from the underlying weighted quiver.


\begin{thebibliography}{99}

\bibitem{DLF_Demonet}
L.~Demonet, \textit{Mutations of group species with potentials and their representations. Applications to cluster algebras}, arXiv:1003.5078.

\bibitem{DLF_DWZ}
H.~Derksen, J.~Weyman, A.~Zelevinsky, \textit{Quivers with potentials and their
representations I: Mutations}, Selecta Mathematica (New Series) \textbf{14} (2008), No. 1, 59--119.

\bibitem{DLF_DR}
V.~Dlab, C.M.~Ringel, \textit{Indecomposable representations
of graphs and algebras}, Memoirs of the American Mathematical Society \textbf{6} (1976), No. 173, 1--57.

\bibitem{DLF_LZ1}
D.~Labardini-Fragoso, A.~Zelevinsky, \textit{Strongly primitive species with potentials I: mutations}, arXiv:1306.3495.



\end{thebibliography}
\end{document}